\documentclass{article}
\usepackage{amsmath,amssymb}
\usepackage[english]{babel}
\usepackage{bm}

\newtheorem{theorem}{Theorem}[section]
\newtheorem{corollary}[theorem]{Corollary}
\newtheorem{lemma}[theorem]{Lemma}

\newtheorem{definition}[theorem]{Definition}

\newtheorem{notation}[subsection]{Notation}

\newtheorem{remark}[theorem]{Remark}

\newenvironment{proof}[1][Proof.]{\begin{trivlist}
\item[\hskip \labelsep {\bfseries #1}]}{\end{trivlist}}

\newcommand{\weight}{\operatorname {weight}}
\newcommand{\Der}{\operatorname {Der}}

\newcommand{\Hom}{\operatorname {Hom}}
\newcommand{\Ext}{\operatorname {Ext}}
\newcommand{\Syz}{\operatorname {Syz}}
\newcommand{\Dc}{\mathcal D}
\newcommand{\cD}{\mathcal D}
\newcommand{\Oc}{\mathcal O}
\newcommand{\cO}{\mathcal O}
\newcommand{\ord}{{\rm ord}}
\newcommand{\gr}{{\rm gr}}
\newcommand{\img}{{\rm Im\,}}
\newcommand{\NN}{{\mathbb N}}
\newcommand{\CC}{{\mathbb C}}
\newcommand{\QQ}{{\mathbb Q}}
\newcommand{\ddelta}{{\widetilde{\delta}}}

\title{A vanishing theorem for a class of logarithmic $\cD$-modules}
\date{January 24, 2007}
\author{F.J. Castro-Jim\'{e}nez, J. Gago-Vargas, M.I. Hartillo-Hermoso and J.M. Ucha}
\begin{document}

\maketitle

\hfil {\em To our beloved friend and colleague Pilar Pis\'{o}n, in
memoriam.}

\vspace{1cm}

\begin{abstract} Let $\cO_X$ (resp. $\cD_X$) be the sheaf of holomorphic functions (resp.
the sheaf of linear differential operators with holomorphic
coefficients) on $X=\CC^n$. Let $D\subset X$ be a locally weakly
quasi-homogeneous free divisor defined by a polynomial $f$. In
this paper we prove that, locally, the annihilating ideal of
$1/f^k$  over $\cD_X$ is generated by linear differential
operators of order 1 (for $k$ big enough). For this purpose we
prove a vanishing theorem for the extension groups of a certain
logarithmic $\cD_X$--module with $\cO_X$. The logarithmic
$\cD_X$--module is naturally associated with $D$ (see Notation
\ref{MlogD}). This result is related to the so called Logarithmic
Comparison Theorem.
\end{abstract}

\section{Introduction}\label{intro_0}

Let us denote by ${\mathcal O}_X$ the sheaf of holomorphic
functions on $X:=\CC^n$ and by $\cD_X$ the sheaf of linear
differential operators with holomorphic coefficients on $X$. A
local section $P$ of $\cD_X$  is a finite sum $P=\sum_{\alpha \in
\NN^n } a_\alpha(x) \partial^\alpha$ where $x=(x_1,\ldots,x_n)$,
$\alpha = (\alpha_1,\ldots,\alpha_n)\in \NN^n$, $a_\alpha(x)$ is a
local section of the sheaf $\cO_X$ and $\partial^\alpha$ stands
for $\partial_1^{\alpha_1} \cdots
\partial_n^{\alpha_n}$ each $\partial_i$ being the partial
derivative with respect to the variable $x_i$. The order of such
an element $P$ is by definition the non negative integer
$\ord(P):=\max \{\vert \alpha \vert := \sum_i \alpha_i \, \vert \,
a_\alpha(x) \neq 0\}$.  For each point $p\in X$ we will write
$\cO_p:=\cO_{X,p}$ and $\cD_p:=\cD_{X,p}$.

Let us fix a point  $p\in X$. Denote by $Der({\mathcal O}_p)$ the
${\mathcal O}_p$-module of ${{\CC}}$-derivations of ${\mathcal
O}_p$. The elements in $Der({\mathcal O}_p)$ are called (germs of)
{\it vector fields} at the point $p$. This yields to the sheaf
$Der(\cO_X)$ of vector fields on $X$. Vector fields are linear
differential operators of order 1.

Let $D$ be a divisor (i.e. a hypersurface) on $X$. Following K.
Saito \cite{Saito}, a (germ of) vector field $\delta \in
Der({\mathcal O}_p)$ is said to be {\it logarithmic} with respect
to $D$ if $\delta(f)=af$ for some $a\in {\mathcal O}_p$, where $f$
is a local (reduced) equation of the germ $(D,p)\subset (X,p)$.
The ${\mathcal O}_p$-module of logarithmic vector fields (or
logarithmic derivations) with respect to $D$ is denoted by
$Der(-\log D)_p$ and it is closed under the bracket product
$[-,-]$. This yields a coherent ${\mathcal O}$-module sheaf
denoted by $Der(-\log D)$, which is a sub-module sheaf of
$Der(\cO_X)$.

\begin{definition}{\rm \cite{Saito}}
The germ of divisor $(D,p)\subset (\CC^n,p)$ is said to be free if
the $\cO_p$-module $Der(-\log D)_p$ of germs of logarithmic vector
field with respect to $D$ is free (and in this case it is
necessarily of rank $n$). If $(D,p)$ is free we also say that $D$
is free at $p$.  A divisor $D\subset \CC^n$ is said to be free if
the germ $(D,p)$ is free for any $p\in D$.
\end{definition}

{\it Saito's criterion} \cite{Saito} says that a divisor $D \equiv
(f=0)$ is free at a point $p \in D$ if and only if there exists a
basis $\{\delta_1,\ldots, \delta_n\}$ of $Der(-\log D)_p$, say
$\delta_i=\sum_j a_{ij}\partial_j$,  whose determinant
$\det((a_{ij}))$ is equal to $u\cdot f$, for some invertible power
series $u \in {\mathcal O}_{p}$ (i.e. such that $u(p) \neq 0$).
Smooth divisors and normal crossing divisors are free. By
\cite{Saito}, any plane curve $D\subset {\CC}^2$ is a free
divisor.

\begin{notation} \label{MlogD} The quotient $$M^{\log D}:= \frac{\cD_X}{\cD_X Der(-\log
D)}$$ plays a fundamental role in what follows. It is a coherent
left $\cD_X$--module and has been introduced in this context in
\cite{Cald-ens}. We are going to consider  later some others
related $\cD_X$-modules.
\end{notation}

Attached to each germ $(D,p)$ we will also consider the following
two left ideals in $\cD_p$. First of all, the annihilating ideal
$$Ann_{\cD_p}(1/f) = \{P\in
\cD_p {\mbox { such that }} P(1/f)=0\}$$ where  $f$ is a local
(reduced)  equation of the germ $(D,p)$, and the ideal
$$Ann_{\cD_p}^{(1)}(1/f) =
\cD_p \{P\in \cD_p {\mbox { such that }}  P(1/f)=0 {\mbox { and }}
\ord(P)=1 \}.$$

An order 1 operator $P\in \cD_p$ such that $P(1/f)=0$ must have
the form $\delta + \frac{\delta(f)}{f}$ where $\delta$ is a
logarithmic derivation in $Der(-\log D)_p$.

\begin{notation} Let us denote by $\widetilde{M}^{\log D}$ the coherent
$\cD_X$-module with stalks $$(\widetilde{M}^{\log D})_p :=
\frac{\cD_p}{Ann_{\cD_p}^{(1)}(1/f)}$$ for $f$ a local reduced
equation of $(D,p)$. Although the previous quotient module depends
on the reduced equation $f$ of the germ $(D,p)$ they are all
isomorphic for different reduced equations.
\end{notation}

The module $\widetilde{M}^{\log D}$ admits in the free Spencer
case (see \cite{Castro-Ucha-moscu}; see also
\cite{Calderon-Narvaez-dual} for an intrinsic treatment of these
objects) a free resolution (called the {\em logarithmic Spencer
resolution})  analogous to the one of $M^{\log D}$ (see
\cite{Cald-ens}):
$${\mathcal
D}_X\otimes_{\cO_X} \wedge^\bullet \widetilde{Der}(-\log D)
\rightarrow \widetilde{M}^{\log D}\rightarrow 0$$  where
$\widetilde{Der}(-\log D)$ denotes the free $\cO_X$-module whose
stalks are $$\widetilde{Der}(-\log D)_p:=\{\delta +
\frac{\delta(f)}{f}\, \vert \, \delta \in Der(-\log D)_p\}.$$

\section{Weak quasi-homogeneity}\label{intro}
In this paper we will consider a weight vector as  an element
$(w_1,\ldots,w_n)\in \QQ^n$ with non negative coordinates and such
that at least one $w_i$ is strictly positive.

A weight vector $w=(w_1,\ldots,w_n)$ defines a filtration on the
ring $\cO=\CC\{x\}=\CC\{x_1,\ldots,x_n\}$ of convergent power
series with complex coefficients.

If $g(x)=\sum_\alpha g_\alpha x^\alpha$ is a non zero element in
$\cO$ we define its weight or its $w$-order as
$\ord_w(g):=\min\{\alpha \cdot w = \sum_i \alpha_i w_i\, \vert \,
g_\alpha \not= 0 \}$. By definition the $w$-order of $0$ is
$+\infty$.

The so called $w$-filtration on $\cO$, which is a decreasing
filtration, is defined by
$$F_\nu=F_\nu (\cO) := \{g\in \cO \, \vert \, \ord_w(g)\geq \nu\}$$ for
all $\nu \in \QQ$. We have $F_\nu = \cO$ for $\nu < 0$.

The associated graded ring is by definition
$$\gr^w(\cO):= \oplus_{\nu \geq 0} \frac{F_\nu}{F_{\nu+1}}.$$

Let us denote by $r=r(w)$ the number of non zero coordinates of
$w$. By assumption $1\leq r \leq n$. The  graded ring $\gr^w(\cO)$
is isomorphic to a polynomial ring in $r$ variables  with
coefficients in a convergent power series ring in $n-r$ variables.
To this end, assume (applying if necessary a permutation of the
components $(x_1,\ldots,x_n)$) that $w_i>0$ for $i=1,\ldots,r$ and
$w_j=0$ for $j=r+1,\ldots,n$. Let us write $x'=(x_1,\ldots,x_r)$,
$x''=(x_{r+1},\ldots,x_n)$ and define $$\CC\{x''\}[x']_\nu =
\left\{\sum_{\beta \in \NN^r} g_\beta(x'')(x')^\beta \in
\CC\{x''\}[x']\, \vert \, \sum_i \beta_i w_i = \nu {\mbox{ if }}
g_\beta \not=0\right\}.
$$

Then the vector space $F_\nu /F_{\nu +1}$ is isomorphic to
$\CC\{x''\}[x']_\nu$. In this way, the weight vector $w$ induces a
graded structure on the ring $\CC\{x''\}[x'] = \oplus_{\nu\geq
0}\CC\{x''\}[x']_\nu$.  The graded rings $\gr^w(\cO)$ and
$\CC\{x''\}[x']$ are then isomorphic as graded rings.

If no confusion arises elements in $\CC\{x''\}[x']_\nu$ are called
weakly quasi-homogeneous (or WQH) power series of weight $\nu$
(with respect to the weight vector $w$). Any non zero element $g$
in $\cO$ can be written in an unique way as a sum $$g=\sum_{\nu
\geq 0} g_{\nu}$$ where each $g_\nu$ is a WQH power series of
weight $\nu$. If the weight $w$ has no zero coordinates (i.e. if
$r=r(w)=n$) then $F_\nu /F_{\nu +1}$ is isomorphic to $\CC[x]_\nu$
the vector space of quasi-homogeneous (QH) polynomials of weight
$\nu$.

If a  WQH power series $f(x)$ has strictly positive weight
$\nu=\ord_w(f)$ then it is WQH of weight 1 with respect to the
weight vector $(w_1/\nu,\ldots,w_n/\nu)$.

\begin{definition} Let $p$ be a point in $\CC^n$. A germ of divisor $(D,p) \subset
(\CC^n,p)$ is said to be  weakly quasi-homogeneous (WQH) if it can
be defined by a WQH germ of convergent  power series around $p$.
If $U\subset \CC^n$ is a non empty open set and $D\subset U$ is a
divisor, we say that $D$ is locally weakly quasi-homogeneous
(LWQH) if for any point $p\in D$ the germ $(D,p)$ is WQH.
\end{definition}

If $r(w)=n$ then weakly quasi-homogeneity is nothing but classical
quasi-homogeneity and locally weakly quasi-homogeneity coincides
with the notion of locally quasi-homogeneity (see \cite{Trans}),
i.e.  every locally quasi-homogeneous (LQH) divisor  is LWQH. The
reciprocal does not hold. For example, the surface defined in
$\CC^3$ by the polynomial $xy(x+y)(xz+y)$ is LWQH but it is not
LQH (see \cite{Comment}).

The $\cO$--module of germs of holomorphic vector fields
$\Der_\CC(\cO)$ is also filtered with respect to the given weight
vector $w$, just by giving to each variable $x_i$ the weight $w_i$
and the weight $-w_i$ to the partial derivative $\partial_i$.  The
$w$--order of a non zero element  $\delta = \sum_i a_i
\partial_i \in \Der_\CC(\cO)$ is then the (possibly negative)
rational number $$\ord_w(\delta) = \min \{\ord_w(a_i) - w_i\,
\vert \, i=1,\ldots,n\}.$$

A vector field $\delta = \sum_i a_i \delta_i = \sum_{i,\alpha}
a_{i,\alpha }x^\alpha \partial_i \in \Der_\CC(\cO)$ is said to be
WQH of weight (or $w$--order) $\mu\in \QQ$ with respect to the
weight vector $w$ if all monomials $a_{i,\alpha} x^\alpha
\partial_i$ in $\delta$ have the same weight $\mu$, i.e. if
$a_{i,\alpha}\not= 0$ then $\alpha \cdot w -w_i = \mu$. Any non
zero vector field $\delta\in \Der_\CC(\cO)$ can be written in a
unique way as a sum $\delta = \sum_{\mu \in \QQ}\delta_\mu$ where
$\delta_\mu$ is the WQH part of $\delta$ of $w$-order $\mu$.

We denote by $\chi=\sum w_i x_i \partial_i$ the Euler vector field
with respect to $w$. It is WQH of weight 0. If $g\in \cO$ is WQH
of weight $\nu$ then $\chi(g)=\nu g$.

\begin{remark} For any WQH vector field  $\delta$ of weight
$\nu$, a straightforward calculation proves that
$[\chi,\delta]=\chi \delta - \delta \chi = \nu\delta$.

\end{remark}
\section{Two basic lemmata}
Let $(D,0)\subset (\CC^n,0)$ be a germ of a WQH free divisor,
defined by some WQH power series $f\in\mathcal{O}$ with respect to
a weight vector $w=(w_1,\ldots,w_n)\in \QQ^n$ and assume that the
weight of $f$ is 1. We recall that $Der(-\log D)_0$ stands for the
$\cO_0$--module of germs of logarithmic derivations with respect
to $(D,0)$.

\begin{lemma}\label{lema:base-especial}
There exists a basis $\{\delta_1,\ldots,\delta_n\}$ of $Der(-\log
D)_0$ such that:
\begin{enumerate}
 \item $\delta_1=\chi$.
 \item Every $\delta_i$ is WQH with respect to  the weight vector
 $w$ and $\delta_i(f)=0$ for $i \ge 2$.
  \item If we write  $\delta_i=\sum a_{ij}\partial_j$ for some $a_{ij}$ in $\cO$, then
 $\det(a_{ij})=f$.
\end{enumerate}
\end{lemma}

\begin{proof} The result being well known, we include a complete proof for the sake
of completeness. First of all, we have $$ Der(-\log
D)_0=\Theta_f\oplus\mathcal{O}\cdot\chi,
$$
where $\Theta_f$ is the $\cO$--module of vector fields
annihilating  $f$. The above decomposition  follows from the
equality   $\delta = (\delta - \frac{\delta(f)}{f}\chi) +
\frac{\delta(f)}{f}\chi$  and  the fact that
$(\delta-\frac{\delta(f)}{f}\chi)(f)=0$, which holds for any
$\delta\in Der(-\log D)_0$.

As $D$ is free, that is, $Der(-\log D)_0$ is free of rank $n$,
then  $\Theta_f$ is free of rank $n-1$.

As  $f$ is WQH  of weight or $w$-order 1, so $f_i:=\frac{\partial
f}{\partial x_i}$ is WQH  of weight or $w$--order $1-w_i$ (we
consider, as usual, the power series 0 to be WQH of order $\nu$
for any $\nu\in \QQ$).

Let us denote by $\Syz_{\mathcal{O}}(f_1,\ldots,f_n)$ the
$\cO$--module of syzygies among $(f_1,\ldots,f_n)$. The
$\mathcal{O}$-modules $\Theta_f$ and
$\Syz_{\mathcal{O}}(f_1,\ldots,f_n)$ are naturally isomorphic.

Let us write $A$ for the graded ring $\CC\{x''\}[x']$ and
$A_\nu=\CC\{x''\}[x']_\nu$ for $\nu\in \QQ$ (see the introduction
for the notations). Let us consider the $A$--module $\Theta_{A,f}$
of vector fields with coefficients in $A$ annihilating $f$.  The
$A$-modules $\Theta_{A,f}$ and $\Syz_{A}(f_1,\ldots,f_n)$ are
naturally isomorphic. By assumption $f,f_1,\ldots,f_n$ are
homogeneous elements in $A$ (more precisely,  they are WQH power
series of $w$-order $1,1-w_1,\ldots,1-w_n$ respectively). Then the
syzygy $A$-module $\Syz_{A}(f_1,\ldots,f_n)$ (resp.
$\Theta_{A,f}$) is finitely generated and it is also graded with
respect to the weight vector $w$. An element $(a_1,\ldots,a_n) \in
\Syz_A(f_1,\ldots,f_n)$ (resp. $\sum_i a_i \partial_i \in
\Theta_{A,f}$) is WQH of $w$-order $\mu\in \QQ$ if and only if it
satisfies the condition $a_i\in A_{\mu+w_i}$ for $i=1,\ldots,n$.
Then  the $A$--modules $\Theta_{A,f}$ and
$\Syz_{A}(f_1,\ldots,f_n)$ are naturally isomorphic as graded
$A$--modules. In fact this graded structure is induced  on
$\Syz_A(f_1,\ldots,f_n)\subset A^n[-w]$ by the shifted graded
structure on $A^n[-w]$ where $(A^n[-w])_\mu = \sum_i A_{\mu+w_i}$
for any $\mu \in \QQ$. As the inclusion $A\subset \cO$ is flat,
the $\cO$--module $\Syz_{\mathcal{O}}(f_1,\ldots,f_n)$ (and so
$\Theta_f$) has a finite system $\{\eta_1,\ldots,\eta_m\}$ of WQH
generators.

Applying Saito's criterion (see \cite{Saito}) to
$\{\chi,\eta_1,\ldots,\eta_m\}$ we can choose among these vectors
fields $n$ elements generating $Der(-\log D)_0$ and the
determinant of its coefficients being  equal to $uf$, with $u$
inversible in $\cO$. Moreover, $\chi$ must be in this generating
system, so we write the system
$\{\delta_1,\delta_2,\ldots,\delta_n\}$ with $\delta_1=\chi$ and
$\weight(\delta_i)=\nu_i$ for some $\nu_i\in \QQ$, $i=2,\ldots,n$.
If we write $\delta_i=\sum a_{ij}\partial_j$ then
$\weight(a_{ij})=w_j+\nu_i$, so the determinant
$\det((a_{ij}))=uf$ is WQH of weight $\nu=\sum w_i+\sum \nu_i$.

If $u=\sum u_{\mu}$ then $u_{\mu}=0$ for all $\mu\neq \nu-1$, that
is, we have   $u=u_{\nu-1}$, so $\nu-1=0$. Changing $\delta_2$ by
$\frac{1}{u} \delta_2$ we obtain the desired basis.

\end{proof}

\begin{lemma}\label{lema:pesos}
Let $(D,0)\subset (\mathbb{C}^n,0)$ be a germ of free divisor as
before and $\{ \chi, \delta_2, \ldots, \delta_n \}$ a basis of
$Der(- \log D)_0$ as in Lemma \ref{lema:base-especial}, with
$\weight(\delta_i) = \nu_i$. Then for all subset $J \subset \{ 2,
\ldots, n \}$
  $$
  1 - \sum_{j \in J} \nu_j > 0.
  $$
\end{lemma}
\begin{proof}
  Each $\delta_i = \sum_{j=1}^n a_{ij} \partial_j$ is a
  WQH vector field of weight $\nu_i$ and because
  $\weight(\partial_j) = -\omega_j$ it follows that
  $\nu_i + \omega_j = \weight(a_{ij})\ge 0$, if $a_{ij}\not= 0$. Let $\Delta$ be the matrix
  whose rows are the weights of the $a_{ij}$:
  $$
  \Delta  = \left( \begin{array}{cccccc} \omega_1 & \ldots & \omega_r &
  0 & \ldots & 0 \\ \nu_2 + \omega_1 & \ldots & \nu_2 + \omega_r &
  \nu_2 & \ldots & \nu_2 \\ \vdots \\ \nu_n + \omega_1 & \ldots &
  \nu_n + \omega_r & \nu_n & \ldots & \nu_n \end{array} \right).
  $$
  Each summand in the determinant of the matrix $(a_{ij})$ is
  WQH  and since $\weight(f)=1$, at least one summand is non zero and has weight $1$.
  So there exists
  some $i \in \{1, \ldots, r\}$ such that
  $$
  1 = \omega_i + \sum_{j=2}^n (\nu_j + \omega_{i(j)}), \mbox{ with
  } {i(j)} \ne i.
  $$
  If $J \subset \{2, \ldots, n \}$ then
  $$
  1- \sum_{j \in J} \nu_j = \omega_i + \sum_{j \not \in J}
  (\nu_j + \omega_{i(j)}) + \sum_{j \in J, i(j) \ne i}
  \omega_{i(j)} > 0.
  $$
\end{proof}

\begin{remark}\label{remark:spencer_matrices}
In Theorem  \ref{propo:ext=0} we will compute some Ext groups of the
$\cD_p$--module  $(\widetilde{M}^{\log D})_p$ for a class of free
divisors $D$ and  $p\in D$. For this purpose we will use the
logarithmic Spencer  resolution of $(\widetilde{M}^{\log D})_p$ (see
\cite{Cald-ens}; see also \cite{Castro-Ucha-moscu}):
$$\Dc_p \bigotimes_{\Oc_p} {\bigwedge}^\bullet
\widetilde{Der}(-\log D)_p \rightarrow (\widetilde{M}^{\log D})_p
\rightarrow 0$$ whose differential is defined as
$$\phi_{\ell} (P \otimes \ddelta_1 \wedge \cdots \wedge \ddelta_{\ell} ) =
\sum_{i=1}^\ell (-1)^{i-1} P\ddelta_i \otimes {\ddelta_1} \wedge
\cdots \widehat{()_i} \cdots \wedge \ddelta_\ell
$$
$$ + \sum_{1\leq i<j \leq \ell } (-1)^{i+j} P \otimes [\ddelta_i, \ddelta_j] \wedge
{\ddelta_1} \wedge \cdots  \widehat{()_i} \cdots \widehat{()_j}
\cdots \wedge \ddelta_\ell$$ where $\widehat{()_i},
\widehat{()_j}$ means that corresponding elements are missing and
$\widetilde{\delta}$ is nothing but $\delta+\frac{\delta(f)}{f}$
for any $\delta \in Der(-\log D)_p$, once a local reduced equation
$f$ of $(D,p)$ has been chosen.

We will take a good basis $\{\delta_1,\ldots,\delta_n\}$ of
$Der(-\log D)_p$ as in Lemma \ref{lema:base-especial} which gives
a corresponding basis $\{\ddelta_1,\ldots,\ddelta_n\}$ of
$\widetilde{Der}(-\log D)_p$.

In addition, it will be useful for handling bases in
${\bigwedge}^\ell Der(-\log D)_p$ (and in ${\bigwedge}^\ell
\widetilde{Der}(-\log D)_p$) to consider a lexicographical
ordering with respect to the indexes of the elements $\delta_i$:
$\delta_{i_1} \wedge \cdots \wedge \delta_{i_\ell}$ precedes
$\delta_{j_1} \wedge \cdots \wedge \delta_{j_\ell}$ if $i_1 = j_1,
\ldots, i_s=j_s$ and $i_{s+1} < j_{s+1}$ for some $s<\ell$.

With respect to these bases we will identify
$$\Dc_p \bigotimes_{\Oc_p} {\bigwedge}^\ell Der(-\log D)_p \, \mbox{ and } \,
\Dc_p \bigotimes_{\Oc_p} {\bigwedge}^\ell \widetilde{Der}(-\log
D)_p$$ with $\Dc_p^{n \choose \ell}$. We will also write $\Dc_p^{n
\choose \ell}$ as a direct sum  $R_{\ell} \oplus S_{\ell}$ where
$$R_{\ell} =  \bigoplus _{{2 \leq j_2 <
\cdots < j_{\ell}\leq n}} \cD_p\cdot \ddelta_ 1 \wedge
\ddelta_{j_2} \wedge \cdots \wedge\ddelta_{j_\ell} \,\, {\mbox
{and}} \,\, S_{\ell} = \bigoplus _{{2 \leq i_1 < \cdots <
i_{\ell}\leq n}} \cD_p\cdot \ddelta_{i_1} \wedge \cdots \wedge
\ddelta_{i_\ell}. $$

With this choice of bases the matrices of the morphisms
$\phi_{\ell}: R_\ell \oplus S_\ell \rightarrow R_{\ell-1} \oplus
S_{\ell-1}$ have a special form:
\begin{itemize} \item The coordinate of $\phi_{\ell}(\ddelta_{1}
\wedge \ddelta_{i_2} \wedge \cdots \wedge \ddelta_{i_\ell})$
corresponding to its ``tail" $\ddelta_{i_2} \wedge \cdots \wedge
\ddelta_{i_\ell}$ is $$\ddelta_1 - \sum_{\scriptsize j\in
\{i_2,\ldots,i_\ell\}} \nu_j.$$ \item The coordinate of
$\phi_\ell(\ddelta_{1} \wedge \ddelta_{i_2} \wedge \cdots \wedge
\ddelta_{i_\ell})$ corresponding to the element $\ddelta_{j_2}
\wedge \cdots \wedge \ddelta_{j_\ell}$ is zero if $1 \notin
\{j_2,\ldots,j_\ell\}$ and $(j_2,\ldots,j_\ell) \not=
(i_2,\ldots,i_\ell)$.
\end{itemize}

So the matrices written by rows of the morphisms have the form $$
\left(
\begin{array}{cc}
A_\ell  & X_\ell \\ B_\ell & C_\ell
\end{array} \right),$$
where $X_\ell$ is a diagonal ${n-1 \choose \ell-1} \times {n-1
\choose \ell-1}$ matrix with elements of the form $\ddelta_1 -
\sum_{\scriptsize j\in \{i_2,\ldots,i_\ell\}} \nu_j$ in its
principal diagonal for all $2\le i_2 < \ldots < i_\ell=n$. For our
purposes we do not require to know the particular shape of
matrices $A_\ell, B_\ell, C_\ell$.

Let us remark that, as proven in Lemma \ref{lema:base-especial},
we have  $\ddelta_1 = \delta_1+1=\chi+1$ and $\ddelta_i=\delta_i$
for $i=2,\ldots,n$.
\end{remark}

\begin{remark}\label{remark:euler_vector}
Given an Euler vector field $\chi = \sum_{i=1}^n \omega_i x_i
\partial_i $ with $\omega_i > 0, i=1,\ldots,r$, $\omega_i = 0,
i=r+1,\ldots, n,$ and $\psi \in \Oc = {\CC}\{x_1,\ldots, x_n\}$,
it is clear that the equation $$(\chi + c)(h) = \psi$$ has a
convergent solution $h$ for any $c >0$. To prove that, we
decompose the given power series $\psi$ as the sum of its WQH
parts (with respect to the weight vector $w$) $\psi = \sum_{\nu
\in {\bf Q}^+} \psi_{\nu}$. We write the unknown power series $h$
as  $h = \sum_{\nu \in {\bf Q}^+} h_{\nu}$. From the equation
$$(\chi+c)(h) = \sum_{\nu \in {\bf Q}^+} (\nu+c) h_{\nu}
= \sum_{\nu \in {\bf Q}^+} \psi_{\nu}$$ we get  $h_{\nu} =
\frac{1}{c + \nu}\psi_{\nu}$ and so $h$ is the unique  convergent
solution. Moreover, $(\chi + c)(h) = 0$ implies $h = 0$ and then
the solution $h$ of the non-homogeneous equation $(\chi +
c)(h)=\psi$ is unique once $\psi $ is fixed. More generally, if
$c_1, \ldots,c_r$ are strictly positive real numbers, the
$\CC$--linear morphism from $\cO^r$ to  $\cO^r$ defined by the
diagonal matrix
$$\left(\begin{array} {ccc} \chi +c_1 & & \\ &\ddots & \\ & & \chi +
c_r\end{array}\right)$$ is an isomorphism. This fact will plays a
crucial role in what follows (see the proof of Theorem
\ref{propo:ext=0}).
\end{remark}

\begin{notation} For any integer $k\ge 0$ we also consider, as in the
Introduction, the coherent $\cD_X$-module $\widetilde{M}^{(k) \log
D }$ with stalks
$$(\widetilde{M}^{(k)\log D})_p :=
\frac{\cD_p}{Ann_{\cD_p}^{(1)}(1/f^k)}$$ where $f$ is a local
reduced equation of the germ $(D,p)$.   For $k=1$ one has
$\widetilde{M}^{(1)\log D} = \widetilde{M}^{\log D}$ (see
Introduction).

We have a natural $\cD_X$--module morphism $\varphi^{k}_{D}:
\widetilde{M}^{(k)\log D} \longrightarrow \cO_X[*D]$ verifying
$\varphi^k_{D,p}(\overline{P})=P(\frac{1}{f^k})$ for any $P\in
\cD_{X,p}$. Once a local reduced equation $f\in \cO_p$ of $(D,p)$
has been chosen, we will write $\widetilde{M}^{(k)\log
f}=\widetilde{M}^{(k)\log D}$.
\end{notation}

\begin{theorem}\label{propo:ext=0}
Given a Spencer free divisor $(D,p)$ defined at $p\in X=\CC^n$ by
a WQH power  series $f$, then
$$\Ext_{\Dc_p}^i((\widetilde{M}^{\log f})_p,\Oc_p)= 0$$ for $i=0,\ldots,
n$.
\end{theorem}

\begin{proof}
We can assume $p=0\in X$. 
The preceding Lemmata have prepared the computations:
\begin{itemize} \item We choose an adapted basis $\{\delta_1 = \chi, \delta_2,\ldots, \delta_n \}$ of
$Der(-\log D)_0$ as in Lemma \ref{lema:base-especial} with
$[\chi+1,\delta_j]=\nu_j \delta_j$ and
$[\delta_i,\delta_j]=\sum_{l=2}^n \nu_{l}^{ij}\delta_l$. \item We
use the logarithmic Spencer resolution of $(\widetilde{M}^{\log
D})_0$ with respect to this basis, so the matrices of the
morphisms in this resolution are like in Remark
\ref{remark:spencer_matrices}.
\item The elements $\chi + 1 - \sum_{j\in \{i_2,\ldots,i_\ell\}}
\nu_j$ in the main diagonal of the upper-right blocks of the
matrices $X_\ell$ (see Remark \ref{remark:spencer_matrices})
verify by Lemma \ref{lema:pesos} that $1 - \sum_{j\in
\{i_2,\ldots,i_p\}} \nu_j
> 0$.
\end{itemize}

To compute the Ext groups, we apply the functor
$\Hom_{\Dc_0}(-,\Oc_0)$ to the corresponding logarithmic Spencer
complex (obtained by truncation of the logarithmic Spencer
resolution of $\widetilde{M} ^{\log f}$)
$$0
\longrightarrow \Dc_0 \stackrel{\phi_n}\longrightarrow \Dc_0^{n
\choose n-1} \stackrel{\phi_{n-1}}\longrightarrow \cdots
\stackrel{\phi_{2}}\longrightarrow\Dc_0^{n \choose 1}
\stackrel{\phi_1}\longrightarrow \Dc_0 \longrightarrow 0,$$
obtaining the complex
$$0
\longrightarrow \Oc_0 \stackrel{\phi_1^*}\longrightarrow \Oc_0^{n
\choose 1} \stackrel{\phi_{2}^*}\longrightarrow \cdots
\stackrel{\phi_{n-1}^*}\longrightarrow\Oc_0^{n \choose n-1}
\stackrel{\phi_{n}^*}\longrightarrow \Oc_0 \longrightarrow 0,$$
using the isomorphism $\Hom_{\Dc_0}(\Dc_0^p,\Oc_0)\simeq \Oc_0^p$.
The morphisms $\phi_i^*$ come from the $\phi_i$ just by applying
the functor $\Hom_{\Dc_0}(-,\Oc_0)$.

Using the notation of Remark \ref{remark:spencer_matrices} we will
write the $\cO_0$--module\\$\bigwedge^\ell\widetilde{\Der}(-\log
f)\simeq \cO_0^{n \choose \ell}$ as a direct sum $G_{\ell} \oplus
H_{\ell}$ where
$$G_{\ell} =  \bigoplus _{{2 \leq j_2 <
\cdots < j_{\ell}\leq n}} \cO_0\cdot \ddelta_ 1 \wedge
\ddelta_{j_2} \wedge \cdots \wedge\ddelta_{j_\ell} \,\, {\mbox
{and}} \,\, H_{\ell} = \bigoplus _{{2 \leq i_1 < \cdots <
i_{\ell}\leq n}} \cO_0\cdot \ddelta_{i_1} \wedge \cdots \wedge
\ddelta_{i_\ell}.
$$ We will identify the isomophic $\cO$--modules $H_{\ell-1}$ and $G_\ell$. The matrix of the $\CC$--linear map
$\phi_\ell^*: G_{\ell-1}\oplus H_{\ell-1}\rightarrow
G_{\ell}\oplus H_{\ell}$ is nothing but
$$\left( \begin{array} {cc} A_\ell & X_\ell \\ B_\ell & C_\ell
\end{array} \right).$$ We have $$\Ext_{\Dc_0}^\ell((\widetilde{M}^{\log D})_0,\Oc_0) =
\frac{\ker (\phi^*_{\ell+1})}{\img (\phi^*_{\ell})}.$$

For $(g,h)\in \ker(\phi_{\ell+1}^*)\subset G_\ell \oplus H_\ell $
we have $A_{\ell+1}g+X_{\ell+1}h=B_{\ell+1}g +C_{\ell+1}h=0$. By
Remark \ref{remark:euler_vector} there is a unique $h'\in
H_{\ell-1}$ such that $X_\ell h' = g$. So,
$\phi_\ell^*(0,h')=(g,h'')$ where  $h''= C_\ell(h')\in H_\ell$. So
$(g,h)-(g,h'')=(0,h-h'')$ belongs to $\ker(\phi_{\ell+1}^*)$. In
particular $X_{\ell+1}(h-h'')=0$ and then, again by Remark
\ref{remark:euler_vector}, $h=h''$. We have proven $(g,h)\in
\img(\phi_\ell^*)$ and then
$\Ext_{\Dc_0}^\ell((\widetilde{M}^{\log D})_0,\Oc_0) =0$.

\end{proof}

\begin{remark}\label{ext=0_for Mklog} The main idea of the proof of Theorem
\ref{propo:ext=0} has its origin in \cite{Ucha-Tesis} (see also
\cite{Castro-Ucha-jsc}). Given an integer $k\geq 1$, if
$\widetilde{M}^{(k)\log f}$ admits a logarithmic Spencer
resolution (analogous to the one of Remark
\ref{remark:spencer_matrices}) then Theorem \ref{propo:ext=0}
holds for $\widetilde{M}^{(k)\log f}$.
\end{remark}

\begin{theorem} \label{ann}
  Let $D\subset X$ be a Spencer free divisor.  We also assume
  that $D$ is LWQH on $X$ (i.e. the germ $(D,p)$ can be defined by a WQH germ of
  holomorphic function in $\cO_p$ for any $p\in D$). Then, as long as $(\widetilde{M}^{(k)\log f})_p$
  admits a logarithmic Spencer resolution, we have
$$
Ann_{\Dc_p}\left(\frac{1}{f^k}\right) =
Ann_{\Dc_p}^{(1)}\left(\frac{1}{f^k}\right)
$$
for any $p\in D$,  any reduced equation $f$ of $(D,p)$ and $k \gg
0$.
\end{theorem}

\begin{proof}
Let us consider the $\cD_X$--module $\cO_X[*D]$ of meromorphic
functions on $X$ with poles along $D$. If $p\in X$ and $f$ is a
local equation of $(D,p)$ we have
$(\cO_X[^*D])_p=\cO_p[\frac{1}{f}]$.

Let $p\in D$ and $-k_0$ be the least integer root of the local
$b$-function $b_{f,p}$ where $f\in \cO_p$ is a local reduced
equation of the germ $(D,p)$. We know that $-n \le -k_0 \le -1$.
Let $k$ be an integer $k \ge k_0$. We have an exact sequence
\begin{equation}\label{eq:res-mtilde}
  0 \to L_{k,p} \to (\widetilde{M}^{(k) \log f})_p \to
  \frac{\Dc_{p}}{Ann_{\Dc_{p}}
  (\frac{1}{f^{k}})} = \Oc_p[\frac{1}{f}] \to 0
  \end{equation}
where $L_{k,p}$ is the kernel of the morphism $\varphi_{D,p}^k$
which is surjective  because  $k \ge k_0$.

By considering the long exact sequence associated to the exact
sequence (\ref{eq:res-mtilde}) we get
$$
  \ldots \rightarrow \Ext^i(\Oc_p[\frac{1}{f}], \Oc_p) \rightarrow
  \Ext^i((\widetilde{M}^{(k) \log f})_p, \Oc_p) \rightarrow
  \Ext^i(L_{k,p},\Oc_p) \rightarrow $$ $$ \rightarrow
  \Ext^{i+1}(\Oc_p[\frac{1}{f}], \Oc_p)
  \rightarrow
  \ldots
$$ where $i\ge 0$ and the $\Ext$ groups have been considered
with respect to the ring $\cD_p$.

Since $p \in D$ the vector space
$\Ext_{\cD_p}^i(\Oc_p[\frac{1}{f}],\Oc_p)$ is equal to $0$  for
$i\ge 0$ (see e.g. \cite[Chap. II, Th. 2.2.4]{Meb-hermann}).

So, from the equality $\Ext_{\cD_p}^i((\widetilde{M}^{(k)\log
f})_p, \Oc_p) = 0$ (see Theorem \ref{propo:ext=0} and Remark
\ref{ext=0_for Mklog}) we get $\Ext_{\cD_p}^i(L_{k,p},\Oc_p) = 0$
for $i\ge 0$. If $p \not \in D$ then $(\widetilde{M}^{(k)\log
f})_p \simeq \Oc_p \simeq \Oc_p[\frac{1}{f}]$ and $L_{k,p}=0$ (see
the exact sequence (\ref{eq:res-mtilde})).

So, we have proved that  $\Ext_{\cD_p}^i(L_{k,p},\Oc_p) = 0$ for
$p\in X$ and $k\geq k_0$.

Since $L_{k,p} = \ker(\varphi^k_D)_p$ and
  $$
  ({\mathcal E}xt_{\Dc_{X}}^i(\ker(\varphi^k_D),\Oc_X))_p \simeq \Ext_{\Dc_{p}}^i(L_{k,p},\Oc_p) = 0
  $$
then the following ${\mathcal{E}xt}$ sheaf vanishes
$$ {\mathcal{E}xt}_{\Dc_X}^i(\ker (\varphi^k_D),\Oc_X) = 0.$$
By \cite[Corollary 11.4.-1]{mebkhout2004} this implies that $\ker
(\varphi^k_D)=0$ for $k\geq k_0$ . This proves the theorem.
\end{proof}

\begin{remark} We do not know if the hypothesis about the LWQH
condition on $f$ is necessary in Theorem \ref{ann}. We notice that
by using  \cite[Corollary 11.4.-1]{mebkhout2004} in the last part
of the proof we do not need to assume the holonomy of
$\widetilde{M}^{(k)\log D}$. Let us also notice that the proof of
Theorem \ref{ann} uses very deep results in $\cD$--module theory:
The Grothendieck Comparison Theorem (as presented in \cite[Chap.
II, Th. 2.2.4]{Meb-hermann}) and the biduality Theorem for
$\cD^{\infty}$--modules (as presented in
\cite[11.4.]{mebkhout2004}).
\end{remark}

\begin{corollary}\label{LWQH-Spencer-free-implies-TCL}
Under the hypotheses of Theorem \ref{ann}, if $-1$ is the least
integer root of the local Bernstein polynomial $b_{f,p}(s)$ for
all $p \in D$,  then the Logarithmic Comparison Theorem (LCT)
holds for $D$.
\end{corollary}
\begin{proof} Let us recall that the divisor $D \subset X$ satisfies the Logarithmic
Comparison Theorem (see \cite{Trans}) if the inclusion of the
logarithmic de Rham  complex $\Omega^\bullet(\log D)$ in the
meromorphic de Rham complex $\Omega^\bullet(*D)$ is a
quasi-isomorphism. Under the hypothesis of the Corollary we have
  $$
  Ann_{\Dc_{p}}(\frac{1}{f}) = Ann_{\Dc_{p}}^{(1)}(\frac{1}{f}),
  $$
for any $p\in D$ and then by \cite[Criterion
3.1]{Castro-Ucha-experimental-04} LCT holds for the divisor $D$.
L. Narv\'{a}ez-Macarro pointed out that \cite[Criterion
3.1]{Castro-Ucha-experimental-04} uses
\cite[Cor.4.2]{Calderon-Narvaez-dual}.
\end{proof}

\begin{remark} If $D$ is a locally quasi-homogenous (LQH) free divisor we do not need
to assume that -1 is the least integer root of the Bernstein-Sato
polynomial $b_{f,p}(s)$ for $p\in D$. To this end, any LQH free
divisor is of Spencer type (see \cite[Theorem 3.2
]{Calderon-Narvaez-moscu}) and from the proof of \cite[Theorem
5.2]{Castro-Ucha-moscu} we deduce, for each  $p\in D$, the
equality $Ann_{\Dc_{p}}(\frac{1}{f}) =
Ann_{\Dc_{p}}^{(1)}(\frac{1}{f})$ where $f$ is a local reduced
equation of the germ $(D,p)$. This implies (see e.g.
\cite[Proposition 1.3]{Torrelli}) that -1 is the least  integer
root of $b_{f,p}(s)$ for any $p\in D$. Then by \cite[Criterion
3.1]{Castro-Ucha-experimental-04} LCT holds for $D$.
\end{remark}

\begin{remark} \label{Luis-Narvaez-Macarro} After 
reading the first version of this paper L. Narv\'{a}ez- Macarro told
us that any LWQH free divisor of Spencer type actually satisfies
LCT. This is more general than Corollary
\ref{LWQH-Spencer-free-implies-TCL} because it is not necessary to
assume the condition about the roots of the $b$-function of $f$.
The sketch of his proof is as follows. By \cite[Th.
2.1.1]{Calderon-Narvaez-int-connections} any Koszul free divisor
$D\subset \CC^n$ satisfies LCT if and only if the canonical
morphism \begin{equation}\label{iso_nar_cal} j_!\CC_U
\longrightarrow \Omega^\bullet_X(\log D)(\cO_X(-D))
\end{equation} is an isomorphism in the derived category of
complexes of sheaves of complex vector spaces. Here $X=\CC^n$, $j:
U = X \setminus D \hookrightarrow X$ is the inclusion and
$\Omega^\bullet_X(\log D)(\cO_X(-D))$ is the tensor product of
$\Omega^\bullet_X(\log D)$ with the invertible $\cO_X$--module
$\cO_X(-D)$. In fact, the argument in the proof of \cite[Th.
2.1.1]{Calderon-Narvaez-int-connections} also shows that the
previous  equivalence is also true for free divisors of Spencer
type, because of the  isomorphism $\Omega_X^\bullet (\log D)
\stackrel{\simeq}{\rightarrow} \Omega_X^\bullet (\log
D)(O_X(-D))^{\vee}$, where ${}^{\vee}$ denotes the Verdier's dual,
which is nothing but the intrinsic version of the duality theorem
in \cite[Th. 4.3 ]{Castro-Ucha-moscu}. Moreover, the
quasi-isomorphism (\ref{iso_nar_cal}) holds if and only if the
complex\\ $\Omega^\bullet_X(\log D)(\cO_X(-D))_p$ is exact at any
point $p\in D$. The last complex is nothing but
$$f\cO_{X,p} \rightarrow f \Omega^1_{X}(\log D)_p \rightarrow
\cdots \rightarrow f \Omega^n_X(\log D)_p$$ (for $f$ a reduced
equation of the germ $(D,p)$) which is a filtered complex using
the weight $w$. By an argument of \cite[Lemma
3.3,6]{Mond-free-and-almost-free} (also used in \cite[Section 2
]{Trans}) this complex is quasi-isomorphic to its subcomplex of
weight 0 which is in fact 0 because the weight of $f$ is 1 and any
logarithmic differential form has a non negative weight.

This proves that any LWQH free divisor of Spencer type satisfies
the Logarithmic Comparison Theorem and we do not need to assume
the condition of Corollary \ref{LWQH-Spencer-free-implies-TCL}
about the $b$-function of $f$.

The free divisor $D$ defined in the space $M_{n,n+1}(\CC)$ of
$n\times (n+1)$ matrices by the vanishing of the product of the
maximal minors \cite{granger-mond-nieto-schulze} is LWQH and as
shown in {\em loc. cit.} it is not LQH. Nevertheless, $D$
satisfies the so called Global Logarithmic Comparison Theorem
\cite{granger-mond-nieto-schulze}. For $n=3$ the divisor is also
of Spencer type and then it satisfies LCT.  It seems that for
$n\geq 4$ the divisor $D$ is also of Spencer type.
\end{remark}

\noindent {\bf Acknowledgment.-} We would like to thank M.
Granger, Z. Mebkhout, D. Mond and L. Narv\'{a}ez-Macarro for their
very helpful suggestions and comments. We also thank L.
Narv\'{a}ez-Macarro for sending us his result sketched in Remark
\ref{Luis-Narvaez-Macarro}.

\end{document}